\documentclass[12pt]{article}

\normalbaselines
\usepackage{amsfonts}

\catcode `\@=11
\def\openone{\leavevmode\hbox{\small1\kern-3.8pt\normalsize1}}%
\@addtoreset{equation}{section}

\def\section{\@startsection {section}{1}{\z@}{-3.5ex plus -1ex minus
     -.2ex}{2.3ex plus .2ex}{\normalsize\bf}}
\def\subsection{\@startsection{subsection}{2}{\z@}{-3.25ex plus -1ex minus
 -.2ex}{1.5ex plus .2ex}{\normalsize\bf}}

\def\thebibliography#1{\section*{References\markboth
  {REFERENCES}{REFERENCES}}\list
  {[\arabic{enumi}]}{\settowidth\labelwidth{[#1]}\leftmargin\labelwidth
  \advance\leftmargin\labelsep
  \usecounter{enumi}}
  \def\newblock{\hskip .11em plus .33em minus -.07em}
  \sloppy
  \sfcode`\.=1000\relax}
 
\catcode `\@=12

\newtheorem{thm}{Theorem}
\newtheorem{rem}{Remark}
\newtheorem{cor}{Corollary}

\newtheorem{pr}{Proposition}
\newtheorem{com}{Comment}

\font\aa=msam10 scaled\magstep1
\font\ab=msbm10 scaled\magstep1
\font\aba=msbm7 
\font\ai=eusb10 scaled \magstep 1
\newcommand{\gz}{\mbox{$\mathbb Z$}}
\newcommand{\gr}{\mbox{$\mathbb R$}}
\newcommand{\gcp}{\mbox{${\mathbb{CP}}$}}
\newcommand{\gc}{\mbox{$\mathbb C$}}
\newcommand{\gcm}{\mbox{\aba C}}
\newcommand{\Gras}{$\mbox{$G_n({\mathbb C}^{m+n})$}$}
\newcommand{\men}{\mbox{$\widetilde {\bf M}$}} 
\newcommand{\got}[1]{{\mbox{${\mathfrak{#1}}$}}}
\newcommand{\mb}[1]{{\mbox{\boldmath{$#1$}}}}
\newcommand{\gh}{\mbox{\ai H}}
\newcommand{\gl}{\mbox{\ai L}}

\newcommand{\gk}{\mbox{\ai K}}
\newcommand{\gph}{\mbox{{\bf P}(\gh )}}
\newcommand{\df}{\mbox{$:=$}}

\newcommand{\Ka}{K\"ahler}
\newcommand{\sfera}{\mbox{${\cal S}(\gh )$}}
\newcommand{\prend}{\aa \symbol{'003}}
\newcommand{\inee}{\ab \symbol{'044}}
\newcommand{\ine}{\mbox{\inee}}
\newcommand{\pre}{\mbox{\prend}}
\newcommand{\oo}{\openone}
\begin{document}
\vspace*{2.5cm}
\noindent
{ \bf COHERENT STATES AND GEOMETRY}\\[1.3cm]
\noindent
\hspace*{1in}
\begin{minipage}{13cm}
Stefan Berceanu  \\[0.3cm]
Institute for Physics and Nuclear Engineering\\
Department of Theoretical Physics\\
PO BOX MG-6, Bucharest-Magurele, Romania\\
E-mail: Berceanu@theor1.ifa.ro; Berceanu@Roifa.Ifa.Ro\\


\end{minipage}

\begin{abstract}
\noindent

The coherent states are viewed as a powerful tool in differential geometry.
It is shown that some objects in differential geometry can be expressed using
quantities which appear in the construction of the coherent states.
The following subjects are discussed via the coherent states:
 the geodesics,
 the conjugate locus and  the cut locus;  the divisors;  the Calabi's
diastasis and its domain of definition;  the Euler-Poincar\'e characteristic
of the manifold, the number of Borel-Morse cells, Kodaira embeding theorem....

\end{abstract}

\section{\hspace{-4mm}.\hspace{2mm} INTRODUCTION}
The coherent states were invented by Erwin Schr\"odinger in the
physical context of the quantum harmonic oscillator
\cite{sch}.
Lately, they were used in Quantum  Optics by
Roy Glauber, E. C. G. Sudarshan   and  John Klauder
 \cite{klauder}. Robert Gilmore \cite{gilmore} and Askold
Perelomov \cite{per} have introduced a group theoretical
setting  which contains Glauber's coherent states
as the particular case of the Heisenberg-Weyl group.
Perelomov's definition of coherent states was globalized by John
Rawnsley \cite{raw}. Perelomov's coherent state vectors ${\mb e}_{Z}$ are points
in  a Hilbert space \gh~ parametrized by the points $Z$ of the so
called coherent state manifold \men , while Rawnsley's coherent state vectors
${\mb e}_{q}$ are indexed by the points $q$ of the so called quantum
line bundle ${\bf M}$ over \men . In the case of homogeneous
Hodge manifolds \men , the quantum line bundle is the pull-back
of the hyperplane bundle $[1_N]$ on the projective space
$\gcp^N$ in  which \men~ is embedded \cite{cgr}.
However, the transition probability,
i.e. the square of the scalar product of the coherent vectors in \gh ,
depends in both cases on the points of the manifold \men .

Actually, a huge literatature in
Physics and Mathematical Physics is devoted to the coherent states. Here I shall
stress the mathematical impact of the coherent states. I shall show that
using the coherent states some
results in differential geometry can be easily recoverd.

More precisely, in this paper I show how the coherent states permit to find:
1) the geodesics;
2) the conjugate locus; 3) the cut locus; 4) the divisors; 5) the Calabi's
diastasis and its domain of definition; 6) the Euler-Poincar\'e characteristic
of the manifold, the number of Borel-Morse cells, Kodaira embeding theorem....
These results are true for very particular manifolds, to be stated
precisely in every case. For the moment, I was not able to produce new
mathematical information via the coherent state aproach. However, I find very
important that the mathematical objects in the list above  can be expressed
using
some quantities which have a physical relevance. Even more, some new remarks
have emerged. For the example, the remark {\bf Polar divisor = Cut locus},
proved on naturally reductive spaces, gives a very simple
characterization of the cut locus in
terms of coherent states. Of course, it will be interesting to see how far
these properties can be extended.

The original part of this paper is extracted from the references \cite{sbcag}
-- \cite{sbs}. The program of investigation of coherent states developed here
was
firstly announced in \cite{ber1}. A short version of this work is
contained
in \cite{ber3}. Reference \cite{ber4} contains an illustration of the methods
claimed in
this paper on the complex Grassmann manifold \Gras .  I  remeber
that in \cite{ber4} I have pointed out some open problems
referring to the conjugate locus on
\Gras . Here I  present only
sketches of the proofs,  because
some of them have been  already published, while
the rest are still in preparation (\cite{project,sbs}).

\section{\hspace{-4mm}.\hspace{2mm}  PRELIMINARIES: THE COHERENT STATES}

 Firstly some notation is introduced at 1.).
Then the coherent states are defined at 2.).

1.) Let $\chi$ be a continuous
representation of the group $K$ on the Hilbert space
\gk~ and let us consider the principal bundle
\begin{equation}
\label{pbundle}
K\stackrel{i}{\rightarrow} G\stackrel{\lambda}{\rightarrow}\men,
\end{equation}
where \men~ is diffeomorphic with $G/K$, $i$ is the inclusion and
$\lambda $ is the natural projection $\lambda (g)=gK$. Let
${\bf M}_{\chi}\df
\men \times_{\chi}\gk$, or simply ${\bf M}\df\men \times_K\gk$, be
the $G$-{\it homogeneous vector bundle\/} \cite{bott}
associated by the
character $\chi$ to the principal $K$-bundle (\ref{pbundle}).
Let  $U\subset\men$~ be open. We introduce the notation
\begin{equation}
\label{local}
(G)^{U}=\{g\in G|\pi (g)\psi_0\in U\},
\end{equation}
where $\pi$ is a representation of $G$ whose restriction to $K$ is
$\chi$ and $\psi_0\in \gk$ corresponds to the base point $o\in \men$.
Then the continuous sections of
${\bf M}_{\chi}$ over $U$ are precisely the continuous maps
$\sigma :U\rightarrow G\times_{\chi}\gk$ of the form
\begin{equation}
\label{sectiuni}
\sigma(\pi (g)\psi_0)=[g,e_{\sigma}(g)],~e_{\sigma}:(G)^U
\rightarrow \gk,
\end{equation}
where $e_{\sigma}$ satisfies the ``functional equation'':
\begin{equation}
\label{fec}
e_{\sigma}(gp)=\chi (p)^{-1}e_{\sigma}(g), g\in(G)^U, p\in K.
\end{equation}

For homogeneous holomorphic line bundles \cite{tirao} ($\gk =\gc$) the
functions in eq. (\ref{fec}) are holomorphic.

2.) Let $ \xi :\gh^{\star}=\gh \setminus \{0\}\rightarrow\gph,
~ \xi (\mb{z})=[\mb{z}]$
be the mapping  which associates to the point $\mb{z}$ in the punctured
Hilbert space the linear
subspace $[\mb{z}]$ generated by $\mb{z}$, where $[\lambda\mb{z}]=[\mb{z}],
\lambda\in\gc^{*}$. The hermitian scalar product $(\cdot,\cdot )$ on \gh~
is linear in the second argument.

 Let us consider the principal bundle (\ref{pbundle}) and let us suppose
the existence of a map $e:G\rightarrow \gh^*$  as in eq. (\ref{sectiuni})
with the property (\ref{fec}) but globally
defined, i.e. on the neighbourhood (\ref{local}) $(G)^{\men}$.
Then $e(G)$ is called {\it  family of coherent vectors\/} \cite{raw}.
If there is a morphism of principal bundles \cite{hus},
i.e. the following diagram is commutative,

\begin{eqnarray}
\label{coherent}
\nonumber G &
~\stackrel{e}{\longrightarrow}~
& \gh^\star \\
\lambda\downarrow &  & \downarrow\xi\\ \label{comdiag}
\nonumber\men & ~\stackrel{\iota}{\longrightarrow}~ & \gph
\end{eqnarray}
then $\iota (\men )$  is called {\it family of coherent states corresponding
to the family of coherent vectors} $e(G)$ \cite{raw}.

We restrict ourselves to the case where the mapping $\iota$  is an embedding
in some projective Hilbert space
\begin{equation}\label{emb}
\iota :
\widetilde {\bf M} \hookrightarrow \gph .
\end{equation}
Also we consider only  {\it homogeneous} manifolds \men~\cite{per}.
In this paper we shall restrict ourselves to 
 {\it k\"ahlerian embeddings} (cf. Ch. 8 in \cite{besse})
  $\iota$, i.e.
\begin{equation}\label{pback}
\omega_{\men}=\iota^*\omega_{\gph},
\end{equation}
where $\omega$ is the fundamental two-form (i.e. closed,
(strongly)
non-degenerate)  of the K\"ahler manifold   and $\iota^*$ is the pull-back
of the mapping $\iota$.

The manifold \men ~
is called
{\it coherent state manifold} and the
$G$-homogeneous line bundle ${\bf M}_{\chi}$ is called {\it coherent vector
manifold} \cite{sbcag}. We shall also suppose that the line bundle
${\bf M}$ is {\it very ample}.

Let now $\tilde{\pi}$ be a projective (in physical literature \cite{bargmann}
``ray'') representation  associated to the unitary irreducible representation
$\pi$ and $\tilde{G}$ the group of transformations which leaves invariant the
transition probabilities in the complex separable Hilbert space $\gh$.
If we use the projection
$\xi '=\xi_{|\sfera}$, i.e. $\xi ' :\sfera\rightarrow\gph,~
\xi '(\psi )=\tilde{\psi}=\{e^{i\phi}\psi|\phi\in\gr\}$, where \sfera~ is the
unit sphere in \gh , then $\tilde{\pi}\circ\xi '=\xi'\circ \pi_{|\gh}$. The
triplet $(\tilde{\pi}, \tilde{G}, \gh )$ is a quantum system with symmetry in
the sense of Wigner and Bargmann \cite{wigner,bargmann}.
Then the manifold $\men\approx G/K$
can be realized as the orbit $\men =\{\tilde{\pi}(g)\tilde{\mb{e}}_0|g\in G\}$,
where $K$ is the stationary group of $\tilde{\mb{e}}_0$ and
$\mb{e}_0\in \gh^{\star}$ is fixed.
For a compact connected simply connected Lie group $G$,  the
existence of the representation $\tilde{\pi}$ implies the existence of the
unitary irreducible representation $\pi$ (cf. the theorem of
Wigner and Bargmann \cite{wigner,bargmann}). This implies the
existence of cross sections $\sigma :\men\rightarrow \sfera$. However, the
(Hopf) principal bundle $\mb{\xi}'=(\sfera ,\xi ', \gph )$ is a $U(1)$-bundle
and
in the construction of coherent vector manifold we need line bundles. But the
principal line bundle $\mb{\xi}'$ is obtained from the (tautological)
line bundle $[-1]=\mb{\xi}=(\gh^\star , \xi ,\gph )$ reducing the group
structure from $\gc^*=GL(1,\gc )$ to $U(1)$.

Here we also stress that the theorem of Wigner and Bargmann  is essentially
\cite{emch} the (first) fundamental theorem of projective geometry
\cite{artin}.

In order to have the physical interpretation of the ``classical system''
obtained by dequantizing the quantum one \cite{cval,sbaa},
we have  restricted ourselves to \Ka~ manifolds \men .
For example, for a compact connected simply connected Lie group $G$,
$\men\approx G/K\approx G^{\gcm}/P$ is a  flag 
 manifold \cite{wolf} and the Borel-Weil
theorem \cite{sbw}
assures the  geometrical realisation of the representation $\pi_j$
and of the representation  space $\gh_j$ if $\mb{e}_0=j$. Here $P$
is a parabolic subgroup of the compexification $G^{\gcm}$ of $G$ and $j$ is
the dominant weight.

The representation $\pi_j$ can be uniquely extended to the group homomorphism
$\pi^*_j:~ G^{\gcm}\rightarrow\pi_j^*(G^{\gcm})$, and respectively, Lie algebra
isomorphism
$\stackrel{.}{\pi}^*_j:\got{g}^{\gcm}
\rightarrow\stackrel{.}{\pi}^*_j
(\got{g}^{\gcm})$ by
\begin{equation}
\pi^*_j(\mbox{e}^Z)=\exp ({\stackrel{.}{\pi}^*_j(Z)}), Z\in\got{g}^{\gcm},
\end{equation}
where $\mbox{e}:\got{g}^{\gcm}\rightarrow G^{\gcm}$ and ${\exp}
: \stackrel{.}{\pi}^*_j\rightarrow\pi^*_j(G)$ are exponential maps, while
$\stackrel{.}{\pi}^*_j(\got{g}^{\gcm})$ is the
complexification of the Lie algebra
$\stackrel{.}{\pi}_j(\got{g})$. We use also the notation $F_{\alpha}=
\stackrel{.}{\pi}^*_j(f_{\alpha})$,
where $\alpha $ is in the set $\Delta$ of the
roots of the Lie algebra $\got{g}$ of $G$ with generators $f_{\alpha}$
of the Cartan-Weyl base of $\got{g}^{\gcm}$ (see also \cite{sbcag}).

Then $\mb{e}_g \df e(g)\df \pi^* (g)\mb{e}_0, g\in G^{\gcm}$ is the
family of coherent vectors, while $\{\tilde{\mb{e}}\}_{g\in G^{\gcm}}$
is  the family of coherent states. The relation
$\mb{e}_g={\mbox{e}}^{i\alpha(g)}\mb{e}_{\lambda (g)}$
defines a  fibre bundle with base \men~ and fibre $U(1)$
\cite{per}. More precisely,
{\it the function}
\begin{equation}
\label{functiune}
\mb{\Upsilon} (g)=(\Upsilon,\mb{e}_g )
\end{equation}
{\it   is holomorphic on $G^{\gcm}$ and defines
holomorphic sections
on the homogeneous holomorphic  line bundle ${\bf M}$ associated
to the principal line bundle $P\rightarrow G^{\gcm}\rightarrow G^{\gcm}/P$ by
the holomorphic character} $\chi$
\begin{equation}
\label{car}
\pi^* (p)e_0=\chi^{-1}(p)e_0,~p\in P,
~\chi (p)={\mbox{\rm e}}^{-i\alpha (p)}.
\end{equation}
Indeed, the function $\mb{\Upsilon} (g)$ verifies
$\mb{\Upsilon} (gp)=\chi^{-1}(p)\mb{\Upsilon} (g), g\in G^{\gcm},~p\in P$,
i.e. eq. (\ref{fec}),
and the corresponding  holomorphic sections are associated via
eq. (\ref{sectiuni}).

Let also  the function
\begin{equation}
\label{u1}
\mb{\Upsilon} '(g)=\mb{\Upsilon} '(gP):=
\displaystyle{\frac{\mb{\Upsilon} (g)}
{(\mb{e}_0 ,\mb{e}_g)}},
\end{equation}
defined on the set
\begin{equation}\label{non-zero}
(\mb{e}_0 ,\mb{e}_g)\not= 0.
\end{equation}
Then
\begin{equation}
\label{uu2}
\mb{\Upsilon}':{\cal V}_0\rightarrow\gc,
~\mb{\Upsilon} '(Z)=(\Upsilon ,\mb{e}_{Z,j}),
\end{equation}
where the  Perelomov's coherent vectors  are
\begin{equation}
{\mb e}_{Z,j}=\exp\sum_{{\varphi}\in\Delta^+_n}(Z_{\varphi}F^+_{\varphi})
j ,~~~~\underline{{\mb e}}_{Z,j}=({\mb e}_{Z,j} ,
{\mb e}_{Z,j})^{-1/2}{\mb e}_{Z,j} ,
\label{z}
\end{equation}
\begin{equation}
{\mb e}_{B,j}=\exp\sum_{{\varphi}\in\Delta^+_n}(B_{\varphi}F^+_{\varphi}-{\bar
B}_{\varphi}F^-_{\varphi}) j , \label{b}
~~~{\mb e}_{B,j}\df  \underline{{\mb e}}_{Z,j}.
\end{equation}
Here $\Delta^+_n$ denotes the positive non-compact roots, $ Z\df (Z_ \varphi )
\in {\gc}^n$  are local
coordinates in the maximal neighbourhood
${\cal V}_0 \subset \widetilde {\bf M} $. Also

\begin{equation}
F^\pm_\varphi=\stackrel{.}{\pi}^{*}(f^\pm_\varphi),
~\varphi\in\Delta^+_n,
\end{equation}
where
\begin{equation}\label{doii}
f^\pm_\varphi=\cases {k^\pm_\varphi=ie_{\pm\varphi},& for $
X_n$,\cr ~~\cr e^\pm_\varphi=e_{\pm\varphi}, &for $X_c$.\cr}
\end{equation}
 $e^\pm_\varphi=e_{\pm\varphi}$ are the part of the Cartan-Weyl base
corresponding to {\got m}. Here
\begin{equation}\label{sume}
{\got g}={\got k}\oplus{\got m}
\end{equation}
is the Cartan
decomposition of the Lie algebra {\got g} of {\bf G} and {\got k} is the Lie
algebra of {\bf K}. The subindex $c$ ($n$) in eq. (\ref{doii})
denotes the compact
(respectively, noncompact) manifold $X_{c,n}\approx \men$.
Then $\got{m}$ is identified with the tangent space at $o$, $\men _o$,
and $\men\approx {\mbox e}^{\got{m}}$.

The homogeneous symmetric spaces  are obtained as  

\begin{equation}
X_{n,c}={\mbox e}^{\sum_{{\varphi}\in\Delta^+_n}(B_{\varphi}f^+_{\varphi}-{\bar
B}_{\varphi}f^-_{\varphi})} \,\cdot~o ,
\end{equation}
where $o=\lambda (e)$ and $e$ is the unit element in $G$.

In eqs. (\ref{z}), (\ref{b})
\begin{equation}\label{action}
F^+_{\varphi} j\not= 0,  ~
F^-_{\varphi} j = 0,~ H_i j=j_i j,\label{act}
\end{equation}
where $ \varphi \in\Delta^+_n$,  $H_i=\pi^*(h_i)$,  $\{h_i\}$ is a
base of the Cartan subalgebra and $ i=1,\ldots,{\rm rank}\,G$.

The system \{$\mb{e}_g$\}, $g\in G^{\gcm}$ is overcomplete
\cite{berezin,per,onofri}
and $(\mb{e}_g,\mb{e}_{g'})$, up to a
factor, is a {\it  reproducing kernel} for the holomorphic vector
bundle $\xi_0:{\bf M}\rightarrow\men$ \cite{zb}.

\section{\hspace{-4mm}.\hspace{2mm}THE  GEODESICS}
 Let us consider again  the orthogonal
decomposition (\ref{sume})
of {\got g} with respect to the $B$-form and let also
${\rm Exp}_p:\widetilde {\bf M}_p\rightarrow \widetilde {\bf M}$ be
the geodesic exponential map.

Let us consider the following two conditions

\parbox{1cm}{\it A)}\parbox[t]{132mm}{$ {\rm Exp} _o
=\lambda \circ {\mbox e} \vert _{\got m}~.  $}

\parbox{1cm}{\it B)}\parbox[t]{132mm}{ On the Lie algebra  {\got g} of $ G$
there exists an $ Ad(G)$-invariant,
symmetric, non-degenerate bilinear form $ B$ such that the restriction
of $B$ to the Lie algebra {\got k}  of $ K$  is likewise non-degenerate.}

We point out that {\it if the homogeneous space} $\widetilde {\bf M}\approx
G/K$ {\it verifies} $B)$, {\it then it also verifies} $A)$ (cf.
Corollary 2.5, Thm. 3.5 and Corollary 3.6 Chapter X in Ref. \cite{kn}).
 Indeed, if
${\got g}={\got k}\oplus {\bf m}$ is the orthogonal
decomposition relative  to the $B$-form on {\got g}, then {\got m} is 
canonically identified with the tangent space at $o,~
 \widetilde {\bf M}_o$. $  B) $ implies a (possibly indefinite)
$G$-invariant metric on $\widetilde {\bf M}$. It follows that $G/K$
is reductive, i.e. $[{\got k},{\got k}]\subset {\got k}$ and $[{\got
k},{\got m}]\subset {\got m}$. If $B)$ is true, then $\widetilde {\bf M}$
is {\it naturally reductive} (see p. 202 in Ref. \cite{kn}) and $A)$ is also
verified. 
The symmetric spaces verify besides the conditions of reductive
spaces, the condition $[{\got m},{\got m}]\subset {\got k}$ and, of course,
$A)$ is  verified too (see Thm. 3.2 Ch. Xl in Ref. \cite{kn}).

Thimm \cite{tim}
 furnishes as another examples of homogeneous spaces verifying $B)$,
 besides the symmetric spaces, the Lie groups with bi-invariant metric
and the normal homogeneous spaces (i.e. $B$ is positive definite).
Kowalski \cite{kow}  studied generalised symmetric spaces still
 verifying condition
$A)$. See also the reference \cite{kowv} for more recent results
on naturally reductive spaces and \cite{mont}.

 Now we remember that in Ref. \cite{sbl}
 we did the following Remark, which is in
 fact E. Cartan's theorem (see e.g. Thm. 3.3 p. 208  in Ref. \cite{helg})
on geodesics on symmetric spaces expressed in the coherent state
setting:

\begin{pr}
 {\it The vector} ${\mb e}_{tB,j}=
\exp \stackrel{.}{\pi}^*_j(tB) j\in {\bf M} , B\in {\got m},$
{\it describes trajectories in} {\bf M} {\it corresponding to the image in the
manifold of coherent states} $\widetilde {\bf M}\hookrightarrow {\gph}$
{\it of geodesics through the identity coset element on the symmetric space}
$X\approx G/K$. {\it The dependence} $Z(t)=Z(tB)$ {\it appearing when
one passes from} eq. (\ref{b}) {\it to} eq. (\ref{z}) {\it describes in
${\cal V}_0$ a geodesic.}
\end{pr}
We shall reformulate  Proposition 1 in a way very useful even for
practical calculations. The proof presented below  implies also Thm. 1.

\begin{pr} \label{pr2}
 {\it For an} $n-$ {\it dimensional manifold} $X\approx G/K$
{\it which has Hermitian symmetric space structure, the parameters}
 $B_{\varphi}$
 {\it in formula} (\ref{b}) {\it of normalised coherent states are normal 
coordinates
in the normal neighbourhood} ${\cal V}_0\approx {\gc}^n${\it around the point}
$Z_{\varphi}=0$ {\it on the manifold} $X$.
\end{pr}
{\it Proof.} The Harish-Chandra embedding theorem can be used (cf. e.g.
Ref. \cite{knapp}; see also Ref. \cite{sbl} for the present
context). This theorem asserts that the map
$M^+\times K^{\gcm}\times M^- \rightarrow G^{\gcm}$ given by
 $(m^+,k,m^-)\rightarrow
m^+km^-$ is a complex analytic diffeomorphism onto an open dense subset
of $G^{\gcm}$ that contains $G_n$. Let ${\got m}^{\pm}$ be the $\pm i$
 eigenspaces
of the complex structure $J$ and  $M^{\pm}$  the (unipotent, Abelian)
subgroups of $G^{\gcm}$
corresponding to ${\got m}^{\pm}$. Then, in particular,
 $b:{\got m}^+ \rightarrow X_c=G^{\gcm}/P,~b(X)= {\mbox e}^{X}P$
is a complex analytic diffeomorphism of ${\got
m}^+$ onto a dense subset of $X_c$ (that contains $X_n$) and the
Remark follows because the requirement $A)$ is fulfilled for the
symmetric spaces.~\pre

\section{\hspace{-4mm}.\hspace{2mm} THE CONJUGATE LOCUS}

Another way to reformulate the Proposition 1 is the following:
\begin{thm}
 {\it Let}  $\widetilde {\bf M}$ {\it be a coherent state 
manifold with Hermitian symmetric space structure,
parametrized in} ${\cal V}_0$ {\it around} $Z=0$ {\it as in eqs.} (\ref{z}),
(\ref{b}). {\it Then the
conjugate locus of the point} $ o $ {\it is obtained vanishing the Jacobian
of the exponential map $Z=Z(B)$ and the corresponding transformations of the
chart from ${\cal V}_0$} .
\end{thm}

{\it Proof.} The proof is contained in  Propsition \ref{pr2}:
{\it the dependence} $Z=Z(B),$
with $ B\in {\got m}^+,$ and $ Z$ 
parametrizing  $\widetilde {\bf M}$, obtained passing 
from eq. (\ref{b}) to (\ref{z}) using   (\ref{action})
(the Baker-Campbell-Hausdorff formulas) \cite{sbl},
{\it expresses} in fact {\it the geodesic exponential} ${\rm Exp}_0: 
 \widetilde {\bf M}_0
\rightarrow  \widetilde {\bf M}$. ~\pre

The situation is very transparent in the case of the complex Grassmann
manifold $X_c=G_n({\gc}^{n+m})=SU(n+m)/S(U(n)\times U(m))$
and his noncompact dual $X_n=SU(n,m)/S(U(n)\times U(m))$. There \cite{sbl} 

\begin{eqnarray}
\nonumber X_{n,c} & = & {\mbox e}^{\left(\matrix{0&B\cr
                            \pm B^*&0\cr}\right)}o=
   \left(\matrix{{\rm co}\sqrt{BB^*}&B{\displaystyle 
{ {\rm si}\sqrt{B^*B}\over \sqrt{B^*B}}}\cr
                                 \pm{\displaystyle {{\rm si} \sqrt{B^*B}\over
                         \sqrt{B^*B}}}B^*&{\rm co}\sqrt{B^*B}\cr}\right)o \\
\nonumber &  & \\
 & = & \left(\matrix{{\oo}&Z\cr 0&{\oo}\cr}\right)
\left(\matrix{({\oo}\mp ZZ^*)^{1/2}&0\cr 0&({\oo}\mp Z^*Z)^{1/2}\cr}\right)
\left(\matrix{{\oo}&0\cr \pm Z^*&{\oo}\cr}\right)o\\
\nonumber &  & \\
\nonumber & = & {\mbox e}^{\left(\matrix{0&Z\cr 0&0\cr}\right)}P ,
\end{eqnarray}
where $B^*$ denotes the hermitian conjugate of the matrix $B$.
co is an abbreviation for the circular cosine cos (resp. the hyperbolic
cosine coh) for $X_c$ (resp. $X_n$) and similarly for si. The $ - ~ (+)$
sign
in the equation above corresponds to the compact (resp. noncompact) $X$.

Here $Z$ and $B$ are $n\times m$ matrices related by the relation 

\begin{equation}
Z=B{{\rm ta} \sqrt{B^*B}\over \sqrt{B^*B}} ,\label{geo}
\end{equation}
and ta is an abbreviation for the hyperbolic tangent tgh (resp. the
circular tangent tg) for $X_n$ (resp. $X_c$). 
The dependence $Z=Z(B)$ describes in fact ${\rm Exp}: 
G_n({\gc}^{n+m})_e \rightarrow G_n({\gc}^{n+m})$ in ${\cal V}_0$.
 Indeed, the equation of geodesics for $X_{c,n}$ is \cite{ber4}

\begin{equation}
\frac{d^2Z}{dt^2}-2\epsilon \frac{dZ}{dt}Z^+
(\oo +\epsilon ZZ^+)^{-1}\frac{dZ}{dt}=0~,\label{65}
\end{equation}
where $\epsilon =1~(-1)$ for $X_c$ (resp. $X_n$). It is easy to see that
 (\ref{geo})
verifies (\ref{65}) with the initial condition $\dot{Z}(0)=B$.

 $Z$ and $B$ in the  eq. (\ref{geo})  of geodesics are in the same time
 the parameters describing
 the coherent states in the paramerization given by eq. (\ref{z}) and
 respectively (\ref{b}).

The following theorem summarizes the known facts about the tangent conjugate
locus and conjugate locus in \Gras~ \cite{wrong,sak,ber4}. The relevant fact
for the  present paper is that {\it the conjugate locus can be calculated using
Theorem 1}.   Below  ${\bf O}^{\perp}$ denotes the orthogonal complement of
the $n$-plane ${\bf O}$ in ${\gc}^N$. We also use the notation

\begin{equation}
V^p_l=\left \{Z\in  G_n({\gc}^{n+m})\vert \dim (Z\cap {\gc}^p)
\geq l\right\} ,
\end{equation}
\begin{equation}
W^p_l=V^p_l-V^p_{l+1}=
\left \{Z\in  G_n({\gc}^{n+m})\vert \dim (Z\cap {\gc}^p)
= l\right\} ,
\end{equation}
\begin{equation}
\omega^p_l=(\underbrace{p-l,\ldots ,p-l}_l, \underbrace{m,\ldots ,m}_{n-l}).
\end{equation}

\begin{equation}
V^p_l=Z(\omega ^p_l);~~ W^p_l=Z'(\omega ^p_l)~ .
\end{equation}

The following  sequences of integers are used
\begin{equation}
\omega=\{ 0\leq\omega (1)\leq \ldots \leq\omega (n)\leq m\},
\end{equation}
\begin{equation}
\sigma (i)=\omega (i)+i,~ i=1,\ldots ,n.
\end{equation}

The Schubert varieties are  defined as
\begin{equation}
Z(\omega )=\left\{ X\in G_n({\gc}^{n+m}) \vert
 \dim (X\cap {\gc}^{\sigma (i)})\geq i\right\} .
\end{equation}

 The "jumps"
sequence  is
\begin{equation}
I_{\omega }=\left \{0=i_0<i_1< \ldots <i_{l-1}<i_{l}=n\right\}~,
\end{equation}

\begin{equation}
\omega (i_h)<\omega (i_{h+1}), \omega (i)=\omega (i_{h-1}), i_{h-1}<i
\leq i_h, h=1,\ldots ,l. 
\end{equation}
 The   generic elements of $Z(\omega )$ are
\begin{equation}
Z'(\omega )=\left\{ X\in  G_n({\gc}^{n+m}) \vert
\dim (X\cap {\gc}^{\sigma (i_h)})= i_h,~  i_h\in I_{\omega}\right\} .
\end{equation}

\begin{thm}
\label{sakth}
The tangent conjugate locus $C_0$ of the point ${\bf O}\in$\Gras~ is given by
\begin{equation}
\label{ura}
C_0=\bigcup_{k,p,q,i}ad\,k(t_iH)~,~i=1,2,3;~1\leq p<q\leq r,
~ k\in K,
\end{equation}
where the vector $H\in{\got a}$ 
   is normalised,
\begin{equation}
H=\sum_{i=1}^r h_iD_{i\,n+i},~h_i\in\gr,~\sum h^2_i=1~.
\label{hh}
\end{equation}
 The parameters $t_i,~i=1,2,3$ in eq. (\ref{ura}) are
\begin{equation}
\begin{array}{l@{\:=\:}c@{\:,\:}l}
t_1  & \displaystyle{\frac{\lambda \pi}{|h_p\pm h_q|}}  &
 ~\mbox{\rm multiplicity}~2;\\[2.ex]
t_2  & \displaystyle{\frac{\lambda \pi}{2|h_p|}} & 
  ~\mbox{\rm multiplicity}~1;\\[2.ex]
t_3  & \displaystyle{\frac{\lambda \pi}{|h_p|}}  & 
 ~\mbox{\rm multiplicity}~2|m-n|; ~\lambda\in \gz^{\star}~.
\end{array}
\label{ttt}
\end{equation}

 The conjugate locus of {\bf O} in \Gras  ~is given by the union
\begin{equation}
\label{reun}
{\bf C}_0={\bf C}^W_0\cup {\bf C}^I_0.
\end{equation}

 The following relations are true

\begin{equation}
\label{ect1}
{\bf C}^I_0= \exp \bigcup_{k,p,q}Ad\,k(t_1H)~,
\end{equation}
\begin{equation}
\label{ect2}
{\bf C}^W_0= \exp \bigcup_{k,p}Ad\,k(t_2H)~,
\end{equation}
i.e. exponentiating the vectors of the type $t_1H$ we get the points
of ${\bf C}^I_0$ for which at least two of the stationary angles with {\bf O}
are equal, while the vectors of the type $t_2H$ are sent to the points of
 ${\bf C}^W_0$ for which at least one of the stationary angles with {\bf O}
is $0$ or $\pi /2$.

The ${\bf C}^W_0$ part of the conjugate locus is given by the disjoint union
\begin{equation}
{\bf C}^W_0=\label{72}
\cases {V^m_1\cup V^n_1,&$ n\leq m,$\cr
            V^m_1\cup V^n_{n-m+1},& $n>m,$ \cr}
\end{equation}
 where
\begin{equation}
V^m_1=\label{74.1}
\cases {{\gcp}^{m-1},& {\mbox {\rm for}} $n=1 ,$\cr
              W^m_1\cup W^m_2\cup \ldots W^m_{r-1}\cup W^m_r,&$1<n ,$\cr}
\end{equation}
\begin{equation}
W^m_r=\cases {G_r({\gc}^{\max (m,n)}),&$n\not= m,$\cr
             {\bf O}^{\perp},&$n=m ,$\cr}
\end{equation}
\begin{equation}
V^n_1=\label{72.2}
\cases {W^n_1\cup \ldots \cup W^n_{r-1}\cup {\bf O},&$1<n\leq m,$\cr
             {\bf O}, &$n=1 ,$ \cr}
\end{equation}
\begin{equation}
V^n_{n-m+1}=W^n_{n-m+1}\cup W^n_{n-m+2}\cup\ldots\cup W^n_{n-1}\cup
{\bf O}~,~ n>m ~.\label{73.3}
\end{equation}
\end{thm}

{\it Sketch of the Proof} (See \cite{ber4}).
 The tangent conjugate locus $C_0$  for \Gras~  in
the case $n\leq m$ was obtained by
 Sakai \cite{sak}. Sakai has observed that Wong's result on the conjugate locus
in the manifold \cite{wrong}
is incomplete, i.e. ${\bf C}^W_0\subset{\bf C}_0$ but
${\bf C}^W_0\ine {\bf C}_0=\exp C_0$. The proof of Sakai consists in
solving the eigenvalue equation $R(X,Y^i)X=e_iY^i$
which appears when solving  the Jacobi equation,
 where the
curvature for the symmetric space $X_c=G_c/K$ at $o$ is simply $R(X,Y)Z=
[[X,Y],Z],~X,Y,Z\in\got{m}_c$. Then $q={\rm Exp}_0tX$ is conjugate to $o$ if
$t=\pi \lambda/\sqrt{e_i},~\lambda\in\gz^{\star}\equiv \gz\setminus \{0\}$.

Above {\got a} is the Cartan subalgebra of the
symmetric pair $(SU(n+m),S(U(n)\times U(m)))$ \cite{helg,sak,ber4}
consisting of vectors of the form (\ref{hh})
where $r$ is the symmetric rank of $X_c$ (and $X_n$) and we use the notation
$D_{ij}=E_{ij}-E_{ji},~i,j=1,\ldots ,N.$
$E_{ij}$ is the matrix with entry $1$ on line $i$ and column $j$ and $0$ 
otherwise.
 The results in the complex Grassmann manifold are obtained farther using the
exponential map given by eq. (\ref{geo}).

The same result on the calculation of the tangent conjugate locus can be
obtained \cite{ber4} using Prop. 3.1 p. 294 in the book of
 Helgason \cite{helg}.  This Proposition asserts that  $H\in
{\got a}$ is conjugate with $o$  iff $\alpha (H)\in i\pi\gz^{\star}$
 for some root $\alpha$  which do not vanish identically on {\got a}. The
 eigenvalues
of the  equation $[H,X]=\lambda X,~\forall H\in {\got a},
~ X\in {\got g}^{\gcm}$ 
lead \cite{ber4} to the values given  in equation (\ref{ura}) for the
 parameters $t_1-t_3$.

 The direct proof \cite{ber4}  in the Grassmann manifold
uses in Theorem 1 the dependence $Z=Z(B)$ furnished by eq.
 (\ref{geo})
which gives the geodesics on 
$G_n({\gc}^{n+m})$  and the Jordan's
 stationary angles between two $n-$planes.

The proof \cite{ber4} is done in four steps. a) Firstly, a diagonalization of
 the $n\times
m$ matrix $Z$ is performed. b) Secondly, the Jacobian of a transformation
of complex dimension one is computed. c) The cut locus is reobtained and his
 contribution to the  conjugate locus is taken into account. d) The
non-zero angles are counted using the following property of the stationary
angles: if the $n'~(n)$-plane $Z'_{n'}$
(resp.  $Z_{n}$) are such that $Z'_{n'}
\cap Z_n=Z"_{n"}$, than $n'-n"$ angles of $Z'_{n'}$ and $Z_n$ are different
from $0$ and $n"$ are 0. ~\pre

\begin{com}
\label{com2}
${\bf C_0}^I$ contains as subset the maximal set of mutually isoclinic
subspaces of the  Grassmann manifold, which are the isoclinic spheres, with
dimension given by the solution of the Hurwitz problem.
\end{com}

{\it Proof.} Wong \cite{wwong} has found out the locus of isoclines in
 $G_n({\gr}^{2n})$,
i.e. the maximal subset $B$   of the Grassmann manifold containing {\bf O}
 consisting of points with
the property that every two $n$-planes   of $B$
have {\it all} the stationary angles equal. Two mutually isoclinic $n$-planes
correspond to the situation where the matrix, which has as eigenvalues
the squares of the stationary angles, is a multiple of
 {\oo}.  The results of Wong were generalized by Wolf \cite{ww1}, who has
considered also the complex and quaternionic Grassmann manifolds.
 The problem of maximal mutually isoclinic subspaces is related with the
Hurwitz problem \cite{hur}. Any maximal set of mutually isoclinic $n$-planes is
 analytically  homeomorphic to a sphere (cf. Thm 8.1 in Wong \cite{wwong} and
 Wolf \cite{ww1}), the dimension of the isoclinic spheres being given by the
solution to the Hurwitz problem.~\pre

\section{\hspace{-4mm}.\hspace{2mm} THE CUT LOCUS }
Let $X$ be complete Riemannian manifold. The point  $q$ is in the
{\it  cut locus} ${\bf CL_p}$ of $p\in X$ if $q$ is the nearest point to $p$ on
the geodesic emanating from $p$ beyond which the geodesic ceases to minimize
his arc length (cf. \cite{kn}, see also Ref. \cite{ber2} for more references).

\begin{rem} $ {\mbox{\rm codim}}_{\gcm} {\bf CL}_p \geq 1$.
\label{rem4}
\end{rem}

We call {\it  polar divisor} of $\mb{e}_0$ the set $
\Sigma _0=\left\{ \mb{e} \in e(G)|(\mb{e}_0, \mb{e})=0\right\}$. This
denomination is inspired after Wu \cite{wu}, who used this term in the case
of the complex Grassmann manifold \Gras .

\begin{thm}\label{cutu}
  Let \men~   be a homogeneous manifold
$\widetilde {\bf M}\approx G/K$.  Suppose that there exists a unitary
irreducible representation $\pi_j$ of $G$ such that in a neighbourhood
$ {\cal V}_0$  around $Z=0$ the coherent states are parametrized as
in eq. (\ref{z}).  Then the manifold \men~ can be represented as
the disjoint union

\begin{equation}
\widetilde {\bf M} ={\cal V}_0\cup \Sigma _0.\label{reu}
\end{equation}

Moreover, if the condition  $B)$ is true, then
\begin{equation}
\Sigma _0={\bf CL}_0. \label{clo}
\end{equation}
\end{thm}

{\it Proof.}
We can take $\psi ~= \psi (Z)~ = {\mb e}_{Z}\in {\bf M}$
 such that
the parameters $Z$ are in ${\gc}^n$ as in formula (\ref{z}). Now, the
second relation (\ref{act}) implies that $(0 , \psi )~=1$ for
 $ \psi \in
\xi^{-1}_0 ({\cal V}_0)$. It follows that the equation 

\begin{equation}
 \cos \theta =0, 
\end{equation}
where

\begin{equation}
\cos\theta ={\vert (0, \psi )\vert\over \Vert 0\Vert ^{1/2}\Vert
\psi \Vert ^{1/2}}=\Vert\psi\Vert ^{-1/2}, \label{unghi}
\end{equation}
does not have solutions for $\psi \in \xi^{-1}_0({\cal V}_0) $,
and the representation (\ref{reu}) follows.

To prove relation (\ref{clo}) if $B)$ is true, use is made of Thm. 7.4 and
the subsequent remark at p. 100 from Ref. \cite{kn}.
The theorem essentially says that any
Riemannian manifold $\widetilde {\bf M}$ is the disjoint union of the
cut locus (closed cell) and the largest open cell of $\widetilde {\bf M}$
on which normal coordinates can be defined. But $Z\in {\gc}^n$ for
points of ${\cal V}_0$ corresponding to the largest normal coordinates 
$B\in {\got m}$, because $B)$ implies $A)$.~\pre

\begin{cor}\label{cul}
 Suppose that $\widetilde {\bf M}$ verifies  $ B)$  and admits the embedding
(\ref{emb}). Let $0, Z\in
\men $.  Then  $Z\in {\bf CL}_0$ iff  the Cayley distance
between the images  $\iota (0), \iota (Z)\in \gph$ is  $\pi /2$
\begin{equation}
d_c(\iota (0),\iota (Z))=\pi /2.
\end{equation}
\end{cor}

Here $d_c$ denotes the the hermitian elliptic Cayley distance on the projective
space
\begin{equation}\label{cd}
d_c([{\mb{\omega}} '],[{\mb{\omega}} ])
=\arccos \frac{\vert({\mb{\omega}} ',{\mb{\omega}})\vert}
{\Vert {\mb{\omega}} '\Vert \Vert {\mb{\omega}} \Vert }~.
\end{equation}
\begin{com}
The cut locus is present everywhere one speaks about coherent states
on symmetric spaces.
\end{com}
This assertion is largely explicated in  reference \cite{sberezin}, where
there are discussed the consequences of theorem \ref{cutu} for the references
\cite{berezin,raw,onofri,moreno,cgr}.

We remember the explicit expression
of the cut locus on the complex projective space and Grassmannian.

\begin{rem}
{\it  On $\gcp^n$,  ${\bf CL}_0=\Sigma_0=H_1=\gcp^{n-1}$}.
\end{rem}
{\it Proof.}
Let the notation
 $ {\cal V}_i=\{z\in \gh^*|z_i\not= 0\}$, ${\cal U}_i=\xi ({\cal V}_i)$,
$H_i=\gph\setminus
{\cal U}_i$. The point
$p_0=[1,0,0,\ldots ]\in \gph$ corresponds to the point
$0$ in the Remark. Then
the solution of the equation $([p_0],[z])=0$ is $[z]=[0,\times ,\times
,\ldots ]= H_1=\gcp^{n-1}\subset\gcp^{n} $ for $\gh ={\gc }^{n}$.\pre

\begin{pr}[Wong \cite{wong}] {\it The cut locus of the point }{\bf O}
 {\it is given by}
\begin{eqnarray}
\nonumber {\bf CL}_0 & = & \Sigma_0  =  V^m_1 = Z(\omega^m_1)= Z(m-1,m,\ldots ,
m)\\ & =  & \left\{ X\subset G_n({\gc}^{n+m})
\vert \dim (X\cap {\bf O}^{\perp})\geq 1\right\} .
\end{eqnarray}

{\it The cut locus in \Gras~ is given by those $n-$planes which have 
at least one of the stationary angles $\pi /2$ with the $n-$plane {\bf O}}.
\end{pr}

{\it Proof.} An immediate proof can be obtained using the results of 
Wu
referring to the polar divisor $\Sigma_0$ on the Grassmann
manifold  (see Ch. l in  Ref. \cite{wu}) and the
theorems  characterising the canonical (universal, det)
bundle on $G_n({\gc}^N)$ (see especially Prop. 3.3 Ch. 7 in
Ref. \cite{hus}), which are particularisations of the representation
in Thm. \ref{cutu}. ~\pre

\section{\hspace{-4mm}.\hspace{2mm} THE DIVISORS}

\begin{pr}
{If \men~ is an homogenous algebraic manifold,
then the polar divisor $\Sigma_0=\iota^*H_1$ is a divisor.}
\end{pr}

{\it Proof.} This fact is a consequence  of
standard properties of divisors \cite{gh} (cf. \cite{polar}). \pre

 The  result is true for any algebraic \Ka~ manifold \cite{sbs}.

Let [~] be the functorial homomorphism  between the group of divisors and
the Picard group of  equivalence class of  holomorphic line bundles \cite{gh}
\begin{equation}
[~]:{\mbox{Div}}(\men )\stackrel{\delta^0}{\rightarrow} H^1(\men ,{\cal O}^*).
\end{equation}

\begin{thm}
{\it Let \men~ be a simply connected Hodge manifold admitting the embedding
(\ref{emb}). Let ${\bf M}=\iota^{*}[1]$ be the
unique, up to equivalence, projectively induced line bundle with a given
admissible connection. Then ${\bf M}=[\Sigma_0]$. Moreover, if the homogeneous
manifold \men~ verifies condition B), then ${\bf M} =[{\bf CL}_0]$.
In particular,
the first relation is true for \Ka ian $C$-spaces, while the second one for
hermitian symmetric spaces\/.}
\end{thm}

{\it Proof}. The main part of the proof is based on the following
theorem of Kodaira and Spencer:  For an algebraic manifold there is a
isomorphism of the group ${\mbox{\rm Cl}}(\men )$
of divisor classes with respect to linear equivalence with the Picard group
${\mbox{\rm Pic}}(\men )$, i.e. for
every complex line bundle ${\bf M}$ over an algebraic manifold \men~ there
exists a divisor $D$ such that $[D]={\bf M}$ (cf. \cite{ks2}).  The
next ingredient is the following theorem due to Kostant:
 Let \men~ be a simply connected Hodge manifold.
 Then, up to equivalence,
there exists a unique line bundle with  a given curvature matrix
$\Theta_{\bf M}$ of the hermitian connection, or,
equivalently, with a given admissible connection (Thm. 2.2.1 in \cite{kost}
p. 135). Farther the theorem \ref{cutu} is used.
The information on \Ka ian $C-$ spaces
is extracted from \cite{wolf,project}. \pre

In the same context, we remember
{\it Bertini's theorem:}
 Let ${\bf M}$ be a projectively  induced line
bundle over an algebraic manifold \men . Then there is a
{\it non-singular} divisor
$D$ of \men~ with ${\bf M}=[D]$. Another formulation of Bertini's theorem is:
 A general hyperplane section $S$ of a connected
non-singular algebraic manifold \men ~in $\gcp^N$ is itself non-singular and for
$n\ge 2$, connected\/ (cf. Zariski \cite{zariski}, Akizuki \cite{akizuki},
Hartshorne \cite{hartshorne}).

\begin{com}{Generally, the divisor $\Sigma_0$ is singular because it doesn't
corresponds to a general section in Bertini's theorem.}
\end{com}

{\it Proof.} We illustrate the assertion on the example
furnished by the Grassmannian $G_2(\gc^4)$.
The Pl\"ucker embedding  is $G_2(\gc^4)\hookrightarrow\gcp^5$.
There are six coordinate neighbourhoods ${\cal V}_1-{\cal V}_6$.
 In ${\cal V}_6$: $p_{12}=0; p_{14}p_{23}-p_{13}p_{24}=0$.
This is a cone over a quadric surface whose vertex is the point
$(0,0,0,0,0,1)$ . The hyperplane $p_{12}=0$ is the
 embedded tangent hyperplane
of $G_2(\gc )^4$
 of the line $x_1=x_2=0$ in $\gcp^5$.  A general hyperplane
section of $G_2(\gc )^4$ is not of the form $p_{12}=0$, since by Bertini's
theorem it has to be smooth. See details in \cite{polar,arrondo}.\pre

\section{\hspace{-4mm}.\hspace{2mm} THE CALABI'S DIASTASIS}


We remember   that the {\it Calabi's diastasis} \cite{calabi}
is expressed through the coherent states (cf. \cite{cgr}) as
\begin{equation}\label{ppp}
D(Z',Z)=-2\log\, \vert (\underline{e}_{Z'} ,\underline{e}_{Z})\vert .
\end{equation}
In proving Proposition \ref{dia}  the following Proposition
\cite{ber1,ber8,ber3} is needed.

\begin{pr}\label{tramp}
Let $\underline{\mb e}_Z$ as in (\ref{z}), where
$ Z$ parametrizes the coherent state
 manifold
in the ${\cal V}_0 \subset \widetilde {\bf M} $ and let us suppose that 
the coherent state manifold admits  the k\"ahlerian embedding
(\ref{emb}).
  Then the angle
\begin{equation} 
\theta := \arccos \vert
(\underline{\mb e}_{Z'} ,\underline{\mb e}_Z)\vert,\label{un1}
\end{equation}
is equal to the Cayley distance on the geodesic joining $\iota (Z'),\iota (Z),$
 where $Z', Z \in {\cal V}_0$,
\begin{equation}
\theta = d_c(\iota (Z'),\iota (Z)).\label{unul1}
\end{equation}
 
More generally, it is true the following relation (Cauchy formula)
 
\begin{equation}\label{cauchy}
(\underline{\mb e}_{Z'},\underline{\mb e}_Z)
=(\iota (Z'),\iota (Z)).
\end{equation}
\end{pr}
{\it Proof.} We discuss firstly the case of compact manifolds.
Let $i:\men \rightarrow \gh$ the mapping $i(Z)={\mb e}_Z$.
Then the
 embedding (\ref{emb})
is realised effectively
  by the formula $\iota =\xi\circ i$, i.e.
\begin{equation}
\label{scut}\iota (Z)=[\underline{\mb e}_Z]~.
\end{equation}

Since the manifold \men~ admits the embedding
(\ref{emb}), the
line bundle ${\bf M}$ is a positive one.
 The following theorem \cite{ss,gh,bott,chern}
 is applied:
 Let ${\bf M}'$ be a holomorphic line bundle  on a compact
complex manifold \men . The following conditions are equivalent:
a) ${\bf M}'$ is positive;
b) for all coherent analytic
 sheaves ${\cal S}$ on \men~ there exists a
 positive integer $m_0(\cal{S})$ such that
 $H^i(\men,{\cal S}\otimes {\bf M}'^m)=0$ for $i>0,~m\geq m_0(\cal{S})$
 (the vanishing theorem of Kodaira);
c) there exists a positive integer $m_0$ such that for all $m\geq m_0$,
there is an embedding  $\iota_{\bf M}:\men \hookrightarrow\gcp^{N-1}$ for
some $N\geq D$ such that ${\bf M}={\bf M'}^m$ is projectively induced, i.e.
${\bf M}=\iota ^*[1]$;
d) \men ~ is a Hodge manifold  (the embedding theorem of Kodaira);
e) the fundamental two-form of \men , the curvature matrix  and the first
 Chern class of ${\bf M}'$ are related by the relations
$ \omega=\frac{\sqrt{-1}}{2}\Theta_{{\bf M}'},
~c_1({\bf M}')=\frac{\omega}{\pi}.$
\newline f) moreover, if \men~ is a k\"ahlerian $C$-space,
then  \men~ is a flag manifold.

But the line bundle ${\bf M}$ is already very ample and
the holomorphic map $\iota_{{\bf M}}:\men\hookrightarrow\gcp^{N-1}$
given by
\begin{equation}
\iota_{{\bf M}}=[s_1(m),\ldots ,s_N(m)]\label{scufund}
\end{equation}
is a holomorphic embedding, where $s_1(m),\ldots ,s_N(m)\in
\Gamma (\men ,{\bf M})$ are global sections.
 Remark that eq. (\ref{scut}) realizes the
embedding (\ref{scufund}).


  The very ample holomorphic line bundle
$ {\bf M} $ of coherent vectors is the
 pull-back $\iota^*$ of
the hyperplane bundle $[1]$ of \gph, the dual bundle of the tautological
line bundle on $\gph^*$, i.e. ${\bf M}=\iota^*[1]$. The analytic line
bundle  ${\bf M}$ is {\it projectively induced} (see p. 139 in Ref.
\cite{hirz}). Then the mapping $i$ preserves the scalar product
\begin{equation}\label{cauchy1}
({\mb e}_{Z'}, {\mb e}_Z)=
(i(Z'),i(Z)),
\end{equation}
which imply
\begin{equation}\label{cauchy2}
(\underline{\mb e}_{Z'},\underline{\mb e}_{Z})=
\frac{(\iota (Z'),\iota (Z))}{\|\iota (Z')\|\|\iota (Z)\|},
\end{equation}
and formula (\ref{cauchy}) follows. Eq. (\ref{unul1}) follows from eq.
(\ref{cauchy2}) and the definition (\ref{cd}) of the Cayley distance.

The noncompact case is treated similarly \cite{kobi}. \pre


\begin{pr}\label{dia}
 The diastasis distance $D(Z',Z)$ between $Z', Z \in {\cal V}_0
 \subset \widetilde {\bf M}$ is related to the geodesic distance
 $\theta = d_c(\iota (Z'),\iota (Z))$,
 where $\iota$ is the embedding (\ref{emb}), by the relation
$$D(Z',Z)= -2\log\, \cos\,\theta .$$
 
 If $\widetilde {\bf M}_n$ is noncompact, 
  $\iota ':\widetilde {\bf M}_n \hookrightarrow \gcp^{N-1,1}=
  SU(N,1)/S(U(N)\times U(1))$,
and $\delta _n (\theta
  _n)$ is the length of the geodesic joining $\iota '(Z'),\iota '(Z)~$
 (resp. $\iota (Z'),\iota (Z))$,
then $$\cos\,\theta _n=
  (\cosh\,\delta _n)^{-1}=e^{-D/2}.$$
\end{pr} 
{\it Proof.} The Proposition is a direct consequence of the relation
(\ref{ppp}) and of the  Proposition \ref{tramp}.

\begin{com}
 The relation  (\ref{clo}) furnishes for manifolds of symmetric type, i.e.
 verifying condition B),
   a geometric description of the
 domain of definition of  Calabi's diastasis: for $z$ fixed, $z'\notin
 {\bf CL}_z$.\end{com}

\section{\hspace{-4mm}.\hspace{2mm} THE EULER-POINCAR\'E CHARACTERISTIC,
THE BOREL-MORSE CELLS, KODAIRA EMBEDDING,...}
  \begin{thm}\label{bigtm}
For flag manifolds
 $\widetilde {\bf M} \approx  G/K $, the following 7 numbers are equal:
 
 1) the maximal number of orthogonal coherent vectors;
 
 2) the number of holomorphic global sections of the holomorphic line bundle
 $\bf M$ over \men, which is supposed to be very ample;
 
 3) the dimension of the fundamental  representation in the Borel-Weil theorem;

 4) the minimal $N$ appearing in the Kodaira embedding theorem,
 $\iota :\widetilde {\bf M} \hookrightarrow  \gcp^{N-1}$;
 
 5) the number of critical points of the energy function $f_H$ attached to a 
 Hamiltonian $H$ linear in the generators of the Cartan algebra of G,
 with unequal coefficients;

 6) the Euler-Poincar\'e characteristic
$\chi (\widetilde {\bf M})=
 [W_G]/[W_H]$,  $[W_G]= card\,W_G$,  where $ W_G$ denotes the Weyl group of
  $G$;

 7) the number of Borel-Morse cells which appear in the CW-complex
 decomposition of $\widetilde {\bf M}$.
\end{thm}

{\it Proof.} The equivalence 1) with 2) follows from the definition of
coherent
states as holomorphic global sections in the $G-$homogenous line bundle $M$.
The Weyl group interchanges the origines in the different coordinate
neighbourhoods of the manifold \men .

 Then
2) is equivalent with 3) due to Borel-Weil theorem
which essentially asserts that \cite{sbw}:
 For every irreducible
representation $\pi_j$ of dominant weight $j$ of the compact connected
semisimple Lie group $G$ corresponds on every homogenous K\"ahlerian space
$G/K\approx G^{\gcm}/P_j$   a complete linear system $|D|$.
The representation space $\gh_j$ of the representation $\pi_j$ is the dual of
$\gl (D)$. The associated line
bundle ${\bf M}'$ is ample iff  the  space $G/K\approx G^{\gcm}/P_j$
is strictly associated to the representation $\tilde{\pi}_j$. Note that
${\bf M}'={\bf M}$, because we have supposed that the line bundle {\bf M} is
already very ample.

 The equivalence  of 1) with 4) follows from the realization of the embedding
 (\ref{emb}) under the form (\ref{scufund}) given by eq. (\ref{scut}).

In order to prove the equivalence of 1) with 5), the energy function
(covariant Berezin symbol \cite{berezin}) $f_H:\men \rightarrow\gr$
\begin{equation}\label{energie}
f_H(Z,{\overline{Z}})=({\underline{\mb e}}_{Z},H{\underline{\mb e}}_{Z})
\end{equation}
is attached to the Hamiltonian
\begin{equation}\label{hamham}
H=\sum_{i=1}^{r}\epsilon_iH_i,~ \epsilon =(\epsilon_1,\ldots ,\epsilon_r)\in\gr.
\end{equation}
The operators $H_i$ are defined in eq. (\ref{action}).
Then it was proved in Reference \cite{sbcag} that the energy function
(\ref{energie}) attached to the Hamiltonian (\ref{hamham}) is a Perfect
Morse function in the extended sense.

The equivalence of 1) with 6) is contained in Theorem 2) in the same reference
\cite{sbcag}. Here we stress only that $\chi (G/K)>0$ if and only if
${\mbox{\rm Rank}}~ G={\mbox{\rm Rank}}~ K$
\cite{hs},
i.e. \men~ is a flag manifold. The same Theorem 2) in \cite{sbcag} implies
the equivalence of 1) and 7).\pre

\begin{com}
The Weil prequantization condition is nothing else that
the condition to have a Kodaira embedding,
i.e. the algebraic manifold \men~ to be a Hodge one.
\end{com}

{\bf Acknowledgments}

The author expresses his thanks to the Organizers
of the Third International Workshop on Differential Geometry and its
Applications for the opportunity to present this talk in Sibiu.
Discussions during the Workshop with Professors Vicente
Cortez, Paul Gauduchon, Tudor Ratiu and Lieven Vanhecke
are kindly acknowledged.

\end{document}